\documentclass[reqno,11pt]{amsart}
\usepackage{hyperref, xy}
\usepackage{amssymb}
\usepackage{fullpage}
\xyoption{all}

\newcommand{\rF}{{\rm{F}}} \newcommand{\bbF}{\mathbb{F}} 
\newcommand{\Spec}{{\rm{Spec}\;}}

\newcommand{\f}{\mathfrak }
  \newcommand{\fgl}{\mathfrak{gl}}
\newcommand{\fg}{\mathfrak{g}}  \newcommand{\fh}{\mathfrak{h}}
\newcommand{\fz}{\mathfrak{z}}

\newcommand{\bbZ}{\mathbb{Z}}
\newcommand{\bbN}{\mathbb{N}}

\newcommand{\mcU}{\mathcal{U}} \newcommand{\mcD}{\mathcal{D}} \newcommand{\mcS}{\mathcal{S}}
\newcommand{\mcO}{\mathcal{O}} \newcommand{\mcV}{\mathcal{V}} \newcommand{\mcZ}{\mathcal{Z}}

\newcommand{\la}{\lambda} \newcommand{\La}{\Lambda}
\renewcommand{\a}{\alpha}
\renewcommand{\b}{\beta}
\newcommand{\g}{\gamma}
\newcommand{\varp}{\overline \omega}

\newtheorem{thm}{Theorem}[section] \newtheorem*{thm*}{Theorem}
\newtheorem*{mainthm}{Main Theorem}
\newtheorem{prop}[thm]{Proposition}
\newtheorem{cor}[thm]{Corollary}
\newtheorem{lem}[thm]{Lemma}

\theoremstyle{definition}
\newtheorem{dfn}[thm]{Definition}

\theoremstyle{remark}
\newtheorem{rem}[thm]{Remark}

\begin{document}
\title[Invariant Differential Operators and FCR factors of $\mcU$]{Invariant Differential Operators and FCR factors of Enveloping algebras}
\author{Ian M. Musson}
\address{Ian M. Musson \\
University of Wisconsin - Milwaukee \\
Department of Mathematical Sciences \\
P.O. Box 0413 \\
Milwaukee WI 53201-0413} \email{musson@uwm.edu}
\thanks{The first author was supported in part by NSF Grant DMS-0099923.}

\author{Jeb F. Willenbring}
\address{
Jeb F. Willenbring \\
University of Wisconsin - Milwaukee \\
Department of Mathematical Sciences \\
P.O. Box 0413 \\
Milwaukee WI 53201-0413} \email{jw@uwm.edu}
\thanks{The second author was supported in part by NSA Grant H98230-05-1-0078.}

\date{\today}

\begin{abstract}
If $\fg$ is a semisimple Lie algebra, we describe the prime factors
of $\mcU(\fg)$ that have enough finite dimensional modules. The
proof depends on some combinatorial facts about the Weyl group which
may be of independent interest. We also determine, which finite
dimensional $\mcU(\fg)$-modules are modules over a given prime
factor. As an application we study finite dimensional modules over
some rings of invariant differential operators arising from Howe
duality.
\end{abstract}
\maketitle

\section{Introduction}

\subsection{}
Let $\rF$ denote an algebraically closed field of characteristic
zero, and let $\mcU = \mcU(\fg)$ be the enveloping algebra of a
reductive Lie algebra $\fg$ over $\rF$.  Starting with the
pioneering work of Soergel (\cite{WS}), there has been some effort
to understand the prime spectrum $X = Spec(\mcU)$ (see \cite{BJ},
\cite{P05}, \cite{P}).  The key open problem is to determine the
order relations between prime ideals.

Intersections of prime ideals are important in relating the
algebraic and topological properties of $X$.  For example, if $P \in
X$, then $P$ is primitive if and only if $P$ is locally closed in
$X$ (see \cite{JD}, 8.5.7).  We can formulate this result as
follows, for $P \in X$, set
\[
    X(P) = \{ Q \in X | P \subset Q \}.
\]
(Throughout this paper we use the symbol $\subset$ to denote proper
inclusion.)  Then exactly one of the following holds:
\begin{enumerate}
\item[(1)] $P$ is primitive,
\item[(2)] $P = \bigcap_{Q \in X(P)} Q$
\end{enumerate}

Now set $Y(P) = \{ Q\in X(P)| \dim \mcU/Q < \infty \}$.  In view of
the above it is natural to ask when $P$ is equal to $\bigcap_{Q \in
Y(P)} Q$. This is closely related to the property that $\mcU/P$ is
an FCR-algebra, (we recall the definition shortly).  In fact, it is
easy to see that if $\fg$ is semisimple then $\mcU/P$ is FCR if and
only if one of the following mutually exclusive conditions hold:
\begin{enumerate}
\item[(3)]  $\dim \mcU/P < \infty$,
\item[(4)]  $P =\bigcap_{Q \in Y(P)} Q$.
\end{enumerate}

\subsection{}

Given an $\rF$-algebra $A$, we say that $A$ \emph{has enough finite
dimensional representations} if the intersection of the annihilators
of the finite dimensional $A$-modules is zero. If in addition all
finite dimensional representations are completely reducible $A$ is
an \emph{FCR} algebra (see \cite{KS}).

Let $\fg$ be a semisimple Lie algebra with enveloping algebra
$\mcU(\fg)$.  In this paper, we describe the prime factors of
$\mcU(\fg)$ that are FCR algebras, in terms of Soergel's theory. If
$A$ is a prime factor of $\mcU(\fg),$ then any finite dimensional
$A$-module is a $\fg$-module, and hence completely reducible.  So
the issue is whether or not $A$ has enough finite dimensional
modules.  This last condition has the following topological
interpretation.  For any algebra $A$, let $$\hat A =
\lim_{\leftarrow} \{ A/I| \dim A/I < \infty \},$$ the profinite
completion of $A$. Then the natural ring homomorphism $A \rightarrow
\hat A$ is injective if and only if $A$ has enough finite
dimensional modules.

\subsection{}
We begin with some basic notation used to state the main theorem.
Suppose that $\fg$ is a reductive Lie algebra.  Let $\fg_{ss} =
[\fg, \fg]$, be the semisimple part of $\fg$.  Fix a Cartan
subalgebra $\fh$ of $\fg$ and set $\fh_{ss} = \fh \cap \fg_{ss}$.
Then, $\fh = \fh_{ss} \oplus \fz$ where $\fz$ is the center of
$\fg$.

Let $R$ denote the root system corresponding to the pair $(\fg,\fh)$
with Weyl group $W$.  Choose a set of positive roots, $R^+$, such
that $R = R^+ \cup -R^+$, and let $B = \{\a_1, \a_2, \cdots, \a_n
\}$ denote the simple roots in $R^+$.

Set $\rho := \frac{1}{2}\sum_{\a \in R^+} \a$, and define the ``dot
action'' of $W$ on $\fh^*$ as $w.\xi := w(\xi + \rho) - \rho$ for
$\xi \in \fh^*$ and $w \in W$. The dot action of $W$ on $\fh^*$
induces an action on $\mcS$. Let $\kappa(\;,\;)$ denote the Killing
form on $\fg_{ss}.$ For $\mu \in \fh^*$, let $h_\mu$ be the element
of $\fh_{ss}$ determined by $\kappa(h_\mu, h) = \mu(h)$ for all $h
\in \fh_{ss}$.  The isomorphism $\fh^*_{ss} \simeq \fh_{ss}$ sending
$\mu$ to $h_\mu$ induces a bilinear form $(\;, \;)$ on $\fh^*_{ss}.$
Recall that this form is positive definite on the real span of the
roots.  For $\a \in R$ let $H_\a = \frac{2 h_\a}{(\a,\a)}$ denote
the coroot to $\a$, and set $\a^\vee = \frac{2 \a}{(\a,\a)}$. Let
$\varp_i \in \fh^*_{ss}$ be the fundamental weights, so that
$(\varp_i,\a^\vee_j) = \delta_{i,j} $. We write $H_i$ (resp. $h_i$)
in place of $H_{\a_i}$ (resp. $h_{\a_i}$).  Thus $\mu(H_i) = (\mu,
\alpha_i^\vee)$.

Each $\la \in \fh^*$ determines a subroot system, $R_\la = \{\a \in
R| (\la, \a^\vee) \in \bbZ \}.$ Let $W_\la$ denote the Weyl group of
$R_\la$. Let $B_\la$ be the unique basis of simple roots in $R^+_\la
:= R^+ \cap R_\la$, and set $W^\la := \left \{ w \in W| w(B_\la)
\subseteq R^+\right \}$. The elements of $W^\la$ are left coset
representatives for $W_\la$ in $W$. Thus, $W = W^\la W_\la$.

\subsection{}
To any prime ideal $\Omega$ of $\mcS = \mcS(\fh)$, the symmetric
algebra on $\fh$, we can associate a prime ideal $I_\Omega$ of
$\mcU$, and a ``tautological highest weight" $\la := \la_\Omega$
(see\cite{WS}). We recall the details in Section \ref{section_FCR}.
Any prime ideal of $\mcU$ has the form $I_\Omega$ for a suitable
prime ideal $\Omega$ in $\mcS$. Set $\mcU_\Omega :=
\mcU/{I_\Omega}$.

For $\mu \in \fh^*$, let $L(\mu)$ be the simple module with highest
weight $\mu$.  Let $$P^+ = \{\mu \in \fh^*| \dim L(\mu) < \infty
\}.$$ We say that $\mu$ is \emph{dominant} if $\mu \in P^+$.  We say
that a prime ideal $\Omega$ is \emph{strongly dominant} if
$\mcV(\Omega) \cap P^+$ is Zariski dense in $\mcV(\Omega)$, where
$\mcV(\Omega)$ denotes the zero set of $\Omega$ in $\fh^*$.  If
$\Omega$ is strongly dominant then $\Omega$ is dominant as defined
in \cite{P}, see Lemma \ref{lem_P}.

\begin{mainthm}  Assume $\fg$ is semisimple and $\mcU = \mcU(\fg)$. \
\begin{enumerate}
\item     If $\mcU_\Omega$ is FCR then for some $w \in W$,
$w.\Omega$ is strongly dominant.
\item
Let $\Omega$ be a strongly dominant ideal in $\mcS$, and $w\in W$.\\
The following conditions are equivalent:
\begin{enumerate}
\item $\mcU_{w.\Omega}$ is FCR

\item $I_{w.\Omega} = I_\Omega$

\item $w \in W^{\la}$ where $\la = \la_\Omega$.
\end{enumerate}
\end{enumerate}
\end{mainthm}

If $\fg$ is reductive we can still give sufficient conditions for
$\mcU$ to have enough finite dimensional representations, see
Proposition \ref{prop_dominant->FCR} and Lemma \ref{lem_c->b}.

This paper is organized as follows. In Section
\ref{sec_weyl_group_combinatorics} we prove some results about root
systems and the Weyl group which may be of independent interest.
Then in Section \ref{section_FCR} we collect some results about
prime and primitive ideals in enveloping algebras. The main theorem
is proved in Section \ref{sec_proof}.  In Sections
\ref{sec_invariant_differential_operators}, \ref{sec_CASE_A}, and
\ref{sec_CASE_BC}, we apply our results to some examples of rings of
invariant differential operators related to Howe duality.

\section{Weyl group combinatorics}\label{sec_weyl_group_combinatorics}

We introduce some notation. 
 Let $B_1$ be a
subset of $B$, $R_1$ be the subroot system of $R$ with simple
roots $B_1$ and let $R_1^+$ be the corresponding set of positive
roots.  Let $W_1$ be the Weyl group of $R_1$ and set
\[\rho^\prime = \frac{1}{2} \sum_{\a \in R_1^+} \a. \]

For $w \in W$ define:
\[
    Q(w) = \left \{ \a \in R^+| w\a \in -R^+ \right \},
\]
and set $\ell(w) = |Q(w)|.$  We have (\cite{GW} Lemma 7.3.6),
\begin{equation}\label{equation_*}
    \rho - w^{-1} \rho = \sum_{\a \in Q(w)} \a.
\end{equation}

\begin{lem}\label{lem_rhoprime}  Suppose $w \in W$ is such that $w(R_1) = R_1$ and
set
\[ T(w) = \{\a \in R^+_1| w \a \in - R_1^+\}. \] Then,
\[\rho^\prime - w^{-1} \rho^\prime = \sum_{\a \in T(w)} \a.\]
\end{lem}
\begin{proof}
Note that $R_1^+$ is a disjoint union $R_1^+ = w(R_1^+ \backslash
T(w)) \cup (-w T(w))$.  The rest of the proof is similar to the
proof of (\ref{equation_*}).
\end{proof}

\begin{prop}\label{prop_positive_general}
Suppose that $w \in W$ and $w(R_1) = R_1$. Then $$( \rho^\prime,
\sum_{\a \in Q(w)} \a ) \geq 0,$$ with equality if and only if
$w(B_1) = B_1$.
\end{prop}

We introduce some notation needed for the proof of this result. If
$Q$ is a finite subset of $\fh^*$ we set $\langle Q \rangle =
\sum_{\a \in Q} \a. $ Let

\[ W_0 = \{ w \in W | w(R_1) = R_1 \}. \]

\begin{lem} \label{2.3} For $w \in W_0$ the following are equivalent,
\begin{enumerate}
\item $w(B_1)=B_1$,
\item $w(B_1) \subseteq R_1^+$,
\item $\ell(w s_\a) > \ell(w)$ for all $\a \in B_1$.
\end{enumerate}
\end{lem}
\begin{proof} Clearly (1) implies (2).  For the reverse implication, suppose
 $\b \in R_1^+$ and write $\b = \sum_{\a_ \in B_1} c_\a \a$.
Define $ht(\b) = \sum_\a c_\a$.  From (2) it follows that $ht (w(\b)
) \geq ht (\b )$ and $w(R^+_1) = R^+_1$.  This implies that
$w^{-1}\g \in B_1$ for $\g \in B_1$.  Hence (1) holds. The
equivalence of (2) and (3) follows from \cite{Hu} Lemma 1.6 a), b).
\end{proof}

\begin{cor} \label{cor} Set $T = \{ w \in W | w(B_1) = B_1 \}$.  Then,
\begin{enumerate}
\item $W_0 = W_1 T$ is the semidirect product of the normal subgroup
$W_1$ by $T$.

\item For $w \in W_1$ and $t \in T$ we have $\ell(wt) = \ell(w) +
\ell(t)$.
\end{enumerate}
\end{cor}
\begin{proof}
(1)    Since $W_1$ is generated by reflections, $ s_\a$ with $\a
\in B_1$ and $ts_\a t^{-1} = s_{t( \alpha )}$ for $t \in T,$
 it follows that $T$
normalizes $W_1$.  From the implication (1) $\implies$ (3) in
Lemma \ref{2.3} and \cite{Hu} Proposition 1.10 (c) we  see that
$T$ is a transversal to $W_1$ in $W_0$. It follows easily that
$W_0 =  W_1 T$ is the semidirect product of $W_1$ by
$T$.\\
(2)  This holds by \cite{Hu} Proposition 1.10 (c).
\end{proof}
Recall that if $w \in W$ and $\a \in B$ such that $\ell(w s_\a) =
\ell(w) + 1$ we have $Q(w s_\a) = s_\a Q(w) \cup \{ \a \}$
(\cite{GW} Corollary 7.3.4).  If $w_1, w_2 \in W$ and $\ell(w_1
w_2) = \ell(w_1) + \ell(w_2)$, it follows by induction on
$\ell(w_2)$ that $Q(w_1 w_2) = w_2^{-1} Q(w_1) \cup Q(w_2),$ a
disjoint union. Now for $t \in T$ and $w \in W_1$ we have
$\ell(wt) = \ell(w) + \ell(t)$ by Corollary \ref{cor}. Since $t
\rho^\prime = \rho^\prime$, it follows that
\[ (\rho^\prime, \langle Q(wt) \rangle ) = (\rho^\prime, \langle Q(w) \rangle) + (\rho^\prime,
    \langle Q(t) \rangle ). \]
Hence Proposition \ref{prop_positive_general} follows from the
following result.

Let $v$ be the longest element of $W_1$.  Then, $v(R_1^+) = - R_1^+$
and since $\ell(v)=|R_1^+|$, we have $v(R^+ \backslash R_1^+)  = R^+
\backslash R_1^+$.  Note that $v^2 = 1$, and if $t \in T$, then
since $t^{-1}vt \in W_1$ and $t^{-1}vt(R_1^+) = -R_1^+$ we have $tv
= vt$.  Suppose $t \in T$, and set
\[ Q^+ = \{ \a \in Q(t) | (\rho^\prime, \a) > 0 \}, \]
\[ Q^- = \{ \a \in Q(t) | (\rho^\prime, \a) < 0 \}. \]

\begin{lem} \
\begin{enumerate}
\item If $w \in W_1$ then $Q(w) = T(w)$. Hence, if $w \neq
1$, then $(\rho^\prime, \langle Q(w) \rangle )>0$.
\item The map $\kappa:\a \mapsto v \a$ is a permutation of $Q(t)$, which interchanges $Q^+$ and $Q^-$.
\item If $t \in
T$, then $(\rho^\prime, \langle Q(t) \rangle )=0$.
\end{enumerate}
\end{lem}
\begin{proof}

\noindent (1) Obviously $T(w) \subseteq Q(w)$, and it follows from
 \cite{Hu} Corollary 1.7 and Proposition
1.10 (b)
 that $\ell(w) = |T(w)|$, giving equality.
 The second statement follows from (\cite{GW} Lemma 2.5.12).\\
  (2)  First observe that if $\a \in Q(t)$ then $\a \notin R_1^+$, since
  $t(R_1^+)=R_1^+$ hence $t \a = -\b$ with $\b \in R^+ \backslash
  R_1^+$.  We must show that $v \a$ is in $R^+$ and $tv \a \in
  -R^+$.  The first statement follows since $v(R^+ \backslash R^+_1)=R^+ \backslash R^+_1$, and the second follows from
  $ tv \a = vt \a = -v \b \in -R^+$.  We have shown that $\kappa$ is
  a permutation of $Q(t)$.  Finally, for $\a \in Q(t)$,
\[ (\rho^\prime, v \a) = (v \rho^\prime, \a) = -(\rho^\prime, \a). \]
The result follows easily from this fact. \\
(3) follows immediately from (2).
\end{proof}

Next we relate the dot action of $W$ on $\mcS$ to the usual
action.

\begin{lem}\label{lem_dot_action}
If $h \in \fh$ and $w \in W$ then,
\[ w.h = wh - \sum_{\a \in Q(w)} h(\a). \]
\end{lem}
\begin{proof}
    For $\la \in \fh^*$ we have,
\begin{eqnarray*}
(w.h)(\la) &=& h(w^{-1}.\la) \\
           &=& h(w^{-1}(\la + \rho)-\rho) \\
           &=& (wh)(\la) + h(w^{-1}\rho-\rho).
\end{eqnarray*}
We now apply equation (\ref{equation_*}).
\end{proof}

\section{FCR factors of Enveloping algebras.}\label{section_FCR}
We recall Soergel's work on prime ideals in the enveloping algebra
$\mcU = \mcU(\fg)$ where $\fg$ is a reductive Lie algebra.  Let
$\mcZ$ denote the center of $\mcU$.

If $R$ is a ring we let Spec($R$) denote the prime spectrum of $R$.
For $\Omega \in \Spec \mcS$, define an ideal  $I_\Omega$ in $\Spec
\mcU$ as follows:  Let $\bbF:=Quot(\mcS/\Omega)$ (the quotient field
of the commutative domain defined by $\Omega$).  We will let
$\fg_\bbF =\bbF \otimes_\rF \fg$ denote the Lie algebra obtained
from $\fg$ by extension of scalars.  We write $\la_\Omega$ for the
$\bbF$-linear map from $\fh_\bbF$ to $\bbF$ whose restriction to
$\fh$ is:
\[
    \la_\Omega : \fh \hookrightarrow \mcS \twoheadrightarrow \mcS/\Omega
    \hookrightarrow \bbF.
\]
Let $L(\la_\Omega)$ be the simple highest weight module for
$\fg_\bbF$ with highest weight $\la_\Omega$. Let $I_\Omega$ be the
annihilator of $L(\la_\Omega)$ in $\mcU(\fg)$ and set
$\phi(\Omega) = I_\Omega$.  The maximal ideals in $\mcS$ have the
form, $M_\mu := \{f \in \mcS| f(\mu) = 0 \}$ for $\mu \in \fh^*$,
and we write $I_\mu$ instead of $I_{M_\mu}$, so that $I_\mu$ is
the annihilator of the simple module $L(\mu)$ with highest weight
$\mu$.  We write $\fh_\bbF^*$ for $\fh^* \otimes \bbF =
(\fh_\bbF)^*.$

For $P \in \Spec \mcU$, let $\psi(P) = P \cap Z$ ($\in \Spec
\mcZ$).  Lastly, $\theta$ denotes the map induced by the
Harish-Chandra isomorphism $\mcZ \stackrel{\sim}{\rightarrow}
\mcS^W \subseteq \mcS$. Then (by \cite{WS}, Section 2.1) we have
the following commutative diagram:
\[
\xymatrix{
     \Spec \mcS \ar[rrd]^\phi \ar[dd]_\theta & &                          \\
                                             & & \Spec \mcU \ar[lld]^\psi \\
     \Spec \mcZ                              & &                          \\
                                             & & }
\]
By \cite{WS} Section 2.2, these maps are continuous in the
Jacobson topology on prime ideals.  In particular these maps
preserve inclusions.

\begin{rem}\label{rem_R}
When $\la = \la_\Omega$ we have:
\[ R_\la = \{ \a \in R| H_\a - n \in \Omega \mbox{ for some } n \in \bbZ \}. \]
\end{rem}
In this situation we define $R_\la^+ = \{ \alpha \in R_\la |
H_\alpha - n \mbox{ for some } n \in \bbN \}$.  In \cite{P}, the
ideal $\Omega$ is called \emph{dominant} if $R^+_\la \subseteq R^+$.
\begin{lem}\label{lem_P}
If $\Omega$ is strongly dominant then $\Omega$ is dominant.
\end{lem}
\begin{proof}
If $\Omega$ is not dominant, then
\[
    H_\alpha + (\rho, \alpha^\vee) - m \in \Omega,
\] for some $\alpha \notin R^+$ and $m \in \bbN$.  If $\beta = -\alpha$ and $n = -m$,
then for $\mu \in \mcV(\Omega)$ we have
\[
    (\mu, \beta^\vee) = n - (\rho,\beta^\vee)< 0,
\] so $\mcV(\Omega) \cap P^+ = \emptyset$, and $\Omega$ is not
strongly dominant.
\end{proof}

Note that the sets $R_\la$, $W_\la$, $B_\la$ depend only on the
coset $\La := \la + P(R) \in \fh^*_\bbF / P(R)$.  For a coset
$\La$, let $B_\La = B_\la$ for any $\la \in \La$.

Set
\[\La^+ = \left \{ \la \in \La | (\la + \rho, \a^\vee) \geq 0, \; \forall \a \in B_\La \right \},
\mbox{and} \] \[\La^{++} = \left \{ \la \in \La | (\la + \rho,
\a^\vee) > 0, \; \forall \a \in B_\La \right \}.\]

\begin{lem}[\cite{Ja} Satz 5.16]\label{lem_Ja 5.16} \
    If $w \in W^\la$ then for all $\mu \in \la + P(R)$ we have $I_\mu = I_{w.\mu}$
\end{lem}

We recall the definition of the $\tau$--invariant for primitive
ideals.  For $w \in W$, set
\[ \tau_\La(w) = \left \{  \a \in B_\La | w.\a \in -R^+ \right \}. \]
Next define $X_\la = \{I_{w.\la}|w \in W\}$, and write $2^{B_\La}$
for the poset of subsets of $B_\La$ ordered by inclusion.  We have:
\begin{thm}[\cite{Ja} Satz 5.7]\label{thm_Ja 5.7}
    Let $\La \in \fh^*/P(R)$ and $\la \in \La^{++}$.  Then there is
    a well defined order reversing map:
    \[ \tau_\La : X_\la \rightarrow 2^{B_\La} \]
    such that $\tau_\La( I_{w.\la} ) = \tau_\La(w)$.
\end{thm}

\begin{lem}\label{lem_density}
Let $D$ be a dense subset of an affine algebraic set $X$.  Let:
\[
    X = X_1 \cup X_2 \cup X_3 \cup \cdots \cup X_r
\]
be the decomposition of $X$ into irreducible components.  Then $D
\cap X_i$ is dense in $X_i$ for each $i$.
\end{lem}
\begin{proof}
Consider $i$ such that $1 \leq i \leq r$.  There exists a function
$p_i \in \mcO(X)$ which vanishes on $X_j$ for all $j \neq i$ and
$p_i$ is not identically zero on $X_i$. Let $f \in \mcO(X_i)$ be a
function that vanishes on $D \cap X_i$.  We will show that $f$ is
identically zero on $X_i$.  We may extend $f$ to a function $g \in
\mcO(X)$. The function $g p_i$ vanishes on $D$ and so is identically
zero on $X$. On $X_i$, $g p_i = f p_i$. This means that $f$ vanishes
on the set $D^\prime = \{x \in X_i| p_i(x) \neq 0\}$, which is dense
in $X_i$. Therefore, $f$ vanishes on $X_i$.
\end{proof}

\section{Proof of the main result.}\label{sec_proof}
Until further notice $\fg$ will be a reductive Lie algebra over $F$.
In general, it is unknown when a prime ideal $I_\Omega$ is contained
in a given primitive ideal $I_\mu$. However for $\mu \in P^+$ the
answer is easy.

We recall the definition of Joseph's characteristic variety
\cite{Jo77}.  Let $p : \mcU \rightarrow \mcU(\fh)$ be the projection
defined by the decomposition
\[
    \mcU = \mcU(\fh) \oplus (\f n^- \mcU + \mcU \f n^+).
\]
For $I$ an ideal of $\mcU$, set $\mathbb V(I) = \mcV(p(I))$.

\begin{lem} Suppose $I, J$ are ideals of $\mcU$ and $\Omega$ is a
prime ideal of $\mcS$.
\begin{enumerate}
\item If $I \subseteq J$ then $\mathbb V(J) \subseteq \mathbb
V(I)$.
\item If $\nu \in \mathbb V(I_\Omega)$ then $I_\Omega \subseteq
I_\nu$.
\item  $\mcV(\Omega) \subseteq \mathbb V(I_\Omega)$.  In particular,
$\mu \in \mathbb V(I_\mu)$.
\item $\mathbb V(I_\Omega) \subseteq W.\mcV(\Omega)$.
\end{enumerate}
\end{lem}
\begin{proof}
The proofs are easily adapted from \cite{Jo77} (Lemma on page 103).
\end{proof}

\begin{thm}\label{thm_ideal-containment}
For $\mu \in P^+$, $I_\Omega \subseteq I_\mu$ if and only if
$W.\mu \cap \mcV(\Omega) \neq \emptyset$.
\end{thm}
\begin{proof} Using the Lemma one can justify the steps below.
\begin{description}
\item[($\Rightarrow$)] If $I_\Omega \subseteq I_\mu$, then $\mu \in \mathbb V(I_\mu) \subseteq \mathbb V(I_\Omega)
\subseteq W.\mcV(\Omega)$. Therefore, $W.\mu \cap \mcV(\Omega) \neq
\emptyset$.
\item[($\Leftarrow$)]
If $w.\mu \in \mcV(\Omega)$ for some $w \in W$, then $w.\mu \in
\mathbb V(I_\Omega)$ since $\mcV(\Omega) \subseteq \mathbb
V(I_\Omega)$. Therefore, since $\mu \in P^+$, $I_\Omega \subseteq
I_{w.\mu} \subseteq I_\mu$.
\end{description} \end{proof}
%

\begin{dfn}\label{def_4.2} Set $\La(\Omega) := \{\mu \in P^+| I_\Omega \subseteq I_\mu
\}$, and for each $w \in W,$ set $$\La(\Omega, w) := \{\mu \in P^+|
w.\mu \in \mcV(\Omega) \}.$$

\bigskip
\noindent Note that, $I_\Omega \subseteq \bigcap_{\mu \in
\La(\Omega)} I_\mu$ and equality holds if and only if $U_\Omega$ is
FCR.  Also, by Theorem \ref{thm_ideal-containment} we have:
\[\La(\Omega) = \bigcup_{w \in W} \La(\Omega,w).\]
\end{dfn}

\begin{prop}\label{prop_dominant->FCR} If $\Omega$ is strongly dominant then $U_\Omega$ has enough finite dimensional modules.
\end{prop}
\begin{proof}
If $P^+ \cap \mcV(\Omega)$ is dense in $\mcV(\Omega)$ then
$I_\Omega = \bigcap_{\mu \in P^+ \cap \mcV(\Omega)} I_\mu$ by
\cite{WS} Proposition 1.  Thus by the remark after Definition
\ref{def_4.2}, it suffices to show that $P^+ \cap \mcV(\Omega)
\subseteq \La(\Omega)$ because this will imply $\bigcap_{\mu \in
\La(\Omega)} I_\mu \subseteq \bigcap_{\mu \in P^+ \cap
\mcV(\Omega)} I_\mu$.  The former inclusion follows since if $\mu
\in P^+ \cap \mcV(\Omega)$ then $\Omega \subseteq M_\mu$, so
$I_\Omega \subseteq I_\mu$, because $\phi$ preserves inclusions.
\end{proof}

\begin{lem}\label{lem_c->b}
If $\Omega$ is a strongly dominant ideal in $\mcS$, and $w \in
W^\la$ where $\la = \la_\Omega$ then $I_{w.\Omega} = I_\Omega$.
\end{lem}
\begin{proof}
If $w \in W^\la$, then $I_\la = I_{w.\la}$ by Lemma \ref{lem_Ja
5.16}. Thus, $I_\Omega = I_{w.\Omega}$ by \cite{WS} Theorem 1 part
(ii).
\end{proof}

For the remainder of this section we assume that $\fg$ is
semisimple.

\begin{prop}\label{prop_Bla_in_B}
 If $\Omega$ is strongly dominant and $\la = \la_\Omega$, then $B_\la \subset B$.
\end{prop}
\begin{proof}
Let $B = \{\a_1, \cdots, \a_n\}$. It suffices to show that if $\a
\in R_\la \cap R^+$ and $\a = \sum_{i=1}^n c_i \a_i$ $(c_i \in \bbZ,
c_i \geq 0)$ then $c_i > 0$ implies $\a_i \in R_\la$.  Now $\a \in
R_\la$ means $H_\a - m^\prime \in \Omega$ for some $m^\prime \in
\bbZ$.  If $\mu \in \mcV(\Omega) \cap P^+$ we have $\mu_j = \mu(H_j)
\in \bbN$ for $1 \leq j \leq n$.  Also, $h_\a = \sum c_j h_j$.  We
can write $H_\a = \sum \frac{r_j}{s} H_j$ where $s, r_j \in \bbN$,
$s
>0$ are such that $\frac{r_j}{s}=c_j \frac{(\a_j, \a_j)}{(\a,\a)}$.
Then, $m^\prime = \mu(H_\a) = \sum \frac{r_j}{s} \mu_j$. Hence, if
$c_i > 0$, then $r_i > 0$ and $0 \leq \mu_i \leq \frac{m^\prime
s}{r_i} = m$.

Suppose $c_i>0$, and for $0 \leq p \leq m$, let $\Pi_{i,p} :=
\mcV(H_i - p)$. The above establishes that $\mcV(\Omega) \cap P^+
\subseteq \cup_{p=0}^m \Pi_{i,p}$.  Since $\Omega$ is strongly
dominant, $\mcV(\Omega) \subseteq \bigcup_{p=0}^m \Pi_{i,p}$, but
$\mcV(\Omega)$ is irreducible so $\mcV(\Omega) \subseteq \Pi_{i,p}$
for some $p$. This means $H_i - p \in \Omega$, so $\a_i \in R_\la$.
\end{proof}

To prove an important special case of the main theorem.
Before proceeding however, observe that for any ideal $\Omega
\subseteq \mcS$ and $w \in W$ we have:
\[ \mcV(w \Omega) = w \mcV(\Omega), \hspace{1cm} \mbox{and} \hspace{1cm} \mcV(w.\Omega) = w.\mcV(\Omega).             \]

Furthermore, if $\Omega^\prime$ is an ideal such that
$\mcV(\Omega^\prime) = \mcV(\Omega) + \rho$, then for any $w \in W$,
$\mcV(w \Omega^\prime) = w.\mcV(\Omega)+\rho$. In addition $H_\a - n
\in \Omega$ implies $H_\a - n - (\rho, \a^\vee) \in \Omega^\prime$.

\begin{thm}\label{thm_both_dominant}
    Suppose $\Omega$ and $w . \Omega$ are prime ideals of $\mcS$
    such that $P^+ \cap \mcV(\Omega)$ and $P^+ \cap \mcV(w.\Omega)$ are
    nonempty.  Then $w \in W^\la$, where $\la := \la_\Omega$ is
    the tautological highest weight corresponding to $\Omega$.
\end{thm}
\begin{proof}
Suppose that $\mu \in P^+\cap\mcV(\Omega)$ and $\nu \in
P^+\cap\mcV(w.\Omega)$.   Assume to the contrary that $w \in W$ is
such that there exists $\a \in B_\La$ such that $w \alpha =
-\beta$ for $\beta \in R^+$. Since $\a \in B_\la$ we have $H_\a -
n \in \Omega$ for some $n \in \bbZ$. Because $\mu \in P^+ \cap
\mcV(\Omega)$, it follows that $\mu(H_\a -n)=0$, and $\mu(H_\a)\in
\bbN$.  This implies that $n \in \bbN$. Similarly $H_\b - m \in
w.\Omega$ for some $m \in \bbN$.

Let $\Omega^\prime$ be the ideal in $\mcS$ defined by the closed
algebraic set $\mcV(\Omega) + \rho$. Then $H_\a-n - (\rho,\a^\vee)
\in \Omega^\prime$, and so $-w(H_\a - n - (\rho, \a^\vee)) =H_\b + n
+ (\rho, \a^\vee) \in w\Omega^\prime$. Since $\nu \in \mcV(w.\Omega)
= w.\mcV(\Omega) = \mcV(w\Omega^\prime) - \rho$, we have, $\nu +
\rho \in \mcV(w\Omega^\prime)$.  Therefore,
\[(\nu+\rho)(H_\beta + n + (\rho, \a^\vee)) = 0, \]
since an element of $w\Omega^\prime$ evaluates as 0 on
$\mcV(w\Omega^\prime)$.  This means that:
\[
    m + (\rho, \b^\vee) + n + (\rho, \a^\vee) = 0,
\]
But this is a contradiction because $(\rho, \gamma^\vee)>0$ for all
$\gamma \in R^+$.
\end{proof}

In the next two lemmas we assume that $\Omega$ is strongly dominant
and $\la = \la_\Omega$.  By Proposition \ref{prop_Bla_in_B}, $B_\la
\subseteq B$.  Thus, we may apply the results of Section
\ref{sec_weyl_group_combinatorics} with $B_1$ replaced by $B_\la$
etc. Let $I$ the subset of $ \{1, 2, \cdots, \ell \}$ such that
$B_\la = \{\a_i | i \in I \}$.

\begin{lem}\label{lem_one}
Assume that $\Omega$ is strongly dominant and $\la = \la_\Omega$.
Then if $w\in W$,
$w(R_\la) = R_\la$ and $w(B_\la) \not\subseteq R_\la^+$ then
\[ w.h_{\rho^\prime} = w h_{\rho^\prime} - c.\]
for some $c>0$.
\end{lem}
\begin{proof}By Lemma \ref{lem_dot_action},
\[
w.h_{\rho^\prime} = w h_{\rho^\prime} - h_{\rho^\prime}(\sum_{\a \in
Q(w)} \a).
\]
Since $\Omega$ is strongly dominant $B_\la \subseteq B$ by
Proposition \ref{prop_Bla_in_B}.  Therefore by Proposition
\ref{prop_positive_general}, $c := h_{\rho^\prime}(\sum_{\a \in
Q(w)} \a) > 0$.
\end{proof}

\begin{lem}\label{lem_wB_in_B}
If $w \in W$ and $w.\Omega = \Omega$ then $w(B_\la) = B_\la$.
\end{lem}
\begin{proof}
Suppose that $w.\Omega = \Omega$.   We proceed in two steps.

Step 1. $w(R_\la) = R_\la$: Indeed, if $\a \in R_\la$, then $H_\a
- m \in \Omega$ for some $m \in \bbZ$.  If $w\a = \b$ then by
Lemma \ref{lem_dot_action}, $\Omega$ contains
\[
    w.(H_\a - m) = H_\b - \sum_{\gamma \in Q(w)}(\g, \a^\vee) -
    m,
\]
and this implies that $\b \in R_\la$.

Step 2: Since the $\a_i$ $(i \in I)$ belong to $R_\la$ there are
non-negative rational numbers $a_i$ ($i \in I$) such that:
\[
    h_i - a_i \in \Omega
\]
similarly, $h_{\rho^\prime} - a \in \Omega$ for some non-negative $a
\in \mathbb{Q}$. By Lemma \ref{lem_rhoprime} we have
\[
    w \rho^\prime - \rho^\prime = \sum_{\a \in T(w)} w\a = - \sum_{i \in I} b_i \a_i
\]
for non-negative integers $b_i$.  Hence,
\begin{equation}\label{equation_one}
    w h_{\rho^\prime} = h_{\rho^\prime} - \sum_{i \in I} b_i h_i.
\end{equation}
Assume for a contradiction that $w(B_\la) \not\subseteq B_\la$.
Then  by Lemma \ref{2.3} $w(B_\la) \not\subseteq R_\la^+,$ so by
Lemma \ref{lem_one},
\begin{equation}\label{equation_two}
w.h_{\rho^\prime} = w h_{\rho^\prime} - c
\end{equation}
for some $c >0$. Since $w.\Omega=\Omega$, it follows from equations
(\ref{equation_one}) and (\ref{equation_two}) that $\Omega$
contains,
\begin{eqnarray*}
w.(h_{\rho^\prime} - a) &=& w h_{\rho^\prime} -a -c \\
                      &=& h_{\rho^\prime} - \sum_{i \in I} b_i h_i -a
                      -c.
\end{eqnarray*}

However, $h_i -a_i \in \Omega$ for $i \in I$, so
\[
    h_{\rho^\prime} - \sum_{i \in I} b_i a_i -a -c \in \Omega.
\]
Since $h_{\rho^\prime}-a \in \Omega$ we deduce that
\[
    \sum_{i \in I} b_i a_i + c \in \Omega.
\]
This is a contradiction since $c > 0$ and $b_i, a_i \geq 0$ for all $i \in I$.
\end{proof}

\begin{thm}\label{thm_Ian's Conjecture 3}
If $U_{\Omega'}$ is FCR there exists a strongly dominant ideal
$\Omega$ in $\mcS$ such that $\Omega' = w.\Omega$ for some $w \in
W^{\la}$ (where $\la := \la_{\Omega}$ is the tautological highest
weight of $\Omega$).
\end{thm}
\begin{proof} Let
\[    X := \overline{\La(\Omega')} = X_1 \cup \cdots \cup X_r.
\] define the decomposition of $X$ into irreducible components.
Set $Y_i := \La(\Omega') \cap X_i$.   Since $U_{\Omega'}$ is FCR,
and $\La(\Omega') = \bigcup Y_i$,
\[I_{\Omega'} = \bigcap_{\mu \in \La(\Omega') }I_\mu =
\bigcap_{i=1}^r \bigcap_{\mu \in Y_i} I_\mu.\]

Because $I_{\Omega'}$ is prime, we have that $I_{\Omega'} =
\bigcap_{\mu \in Y_i} I_\mu$, for some $i$ .  Let $\Omega$ denote
the ideal of elements vanishing on $X_i$. Note that $Y_i$ is dense
in $X_i$ for each $i$ since $\La(\Omega')$ is dense in $X$ (see
Lemma \ref{lem_density}). Therefore $I_{\Omega'} = I_{\Omega}$ by
\cite{WS} Proposition 1. This implies that $\theta(I_{\Omega'}) =
\theta(I_{\Omega})$, and therefore by \cite{WS} Theorem 1 part (i)
we have $\Omega' = w.\Omega$ for some $w \in W$.

Because $\La(\Omega) \cap X_i$ is dense in $X_i$, and $\La(\Omega)
\subseteq P^+,$ it follows that $P^+ \cap \La(\Omega)$ is also dense
in $X_i$. Therefore $\Omega$ is strongly dominant. Since $I_\Omega =
I_{ w.\Omega},$ the proof of \cite{WS} Theorem 1 part (ii) shows
that there exists $w_1 \in W$ such that $w\cdot \Omega = w_1.\Omega$
and $ I_{\lambda} = I_{w_1\cdot \lambda}.$ It follows from Lemma
\ref{lem_wB_in_B} that $w_1^{-1} w(B_\la) = B_\la.$ Hence to prove
the result it is enough to show that
 $w_1 \in W^{\la}.$

Write $w_1 = uv$ with $u \in W^{\la}$, $v \in W_{\la}$.  By Lemma
\ref{lem_Ja 5.16} we have $I_{w_1.{\la}} = I_{v.{\la}}$. Hence
$I_{\la} = I_{v.{\la}}$.  By Theorem \ref{thm_Ja 5.7}, $ \tau_\La(v)
= \tau_\La(1).$  This implies that $v=1,$ since $\tau_\La(1)
=\emptyset$ and if $v \neq 1$ and $v=u s_\alpha$ is a reduced
decomposition in $W_\la$, then $\alpha \in \tau_\Lambda(v)$.
Therefore, $w_1 = u \in W^{\la}$.
\end{proof}

\begin{lem}\label{lem_new}
Suppose $\Omega, \Omega_1 \in \Spec \mcS$ and let $\la, \la_1$
respectively be the tautological highest weights.
\begin{enumerate}
\item If $\Omega_1 = v.\Omega$ with $v \in W^\la$ then $B_{\la_1} = v B_\la$.
\item If in addition, $u \in W^{\la_1}$ then $uv \in W^\la$.
\end{enumerate}
\end{lem}
\begin{proof}\
(1) From remark \ref{rem_R} we see that
\[ v(R_\la) = \{ \a \in R| v^{-1} H_\a - m \in \Omega \mbox{ for some } m \in \bbZ \} \]
and one can check that $v(R_\la) = R_{\la_1}$.  Since $v \in W^\la$,
\[v B_\la \subseteq R^+ \cap R_{\la_1} = R_{\la_1}^+. \]
This fact implies that $v R_\la^+ \subseteq R_{\la_1}^+$,
and therefore $R_{\la_1} \cap \bbN \, v B_\la  = R_{\la_1}^+$.
Since there is a unique choice of simple roots in $R_{\la_1}^+$,
the result follows.

(2) If $\a \in B_\la$, then $v\a \in B_{\la_1}$ by part (1). Hence
$uv \a \in R^+$ for all such $\a$, so $uv \in W^\la$.
\end{proof}

\begin{proof}[Proof of the Main Theorem]\
(1) this follows from Theorem \ref{thm_Ian's Conjecture 3}. \\
For (2) assume $\Omega$ is strongly dominant, and $w \in W$. \\
(b) $\implies$ (a): If $I_{w.\Omega} = I_\Omega$, then
$U_{w.\Omega}$ is FCR by Proposition \ref{prop_dominant->FCR}. \\
(c) $\implies$ (b): This follows from Lemma \ref{lem_c->b}. \\
(a) $\implies$ (c): Suppose that $\mcU_{w.\Omega}$ is FCR.  By
Theorem \ref{thm_Ian's Conjecture 3}, there is a strongly dominant
prime $\Omega_1$, and $u \in W^{\la_1}$ such that $w.\Omega =
u.\Omega_1$, where $\la_1$ is the tautological highest weight
corresponding to $\Omega_1.$  If $v = u^{-1}w$, then $\Omega$ and
$v.\Omega$ are strongly dominant, so by Theorem
\ref{thm_both_dominant}, $v \in W^\la$. By Lemma \ref{lem_new} we
have $w = uv \in W^\la$.
\end{proof}

\section{Invariant differential operators}\label{sec_invariant_differential_operators}
An interesting class of examples arises from dual pairs (\cite{GW}
Section 4.5, \cite{RH}, \cite{LS}). The application of our method is
fairly routine, so we give only the most interesting examples rather
than conduct an exhaustive study.

Let $M_{p,k}$ denote the set of $p \times k$ matrices over $F$, and
consider cases $A$,$B$,$C$ as follows. \vspace{1.0cm}
\begin{center}
 \begin{tabular}{|c|c|c|c|c|} \hline
   {\bf Case} & $\mathbf{K}$
& $\mathbf{V}$   & {\bf Action of} $K$ \bf{on} $\mathbf{V}$ &
$\mathbf{\mathfrak{g}}$ \\ \hline
  $A$  & $GL_k(F)$ &$M_{p,k}\times M_{k,q}$&$(g,(a,b))
                                 \mapsto(ag^{-1},gb)$   & $g\ell_{p+q}$
                                 \\ \hline
  $B$  & $O_k(F)$ & $M_{k,n}$ & $(g,a) \mapsto
                                                    ga$ & $sp_{2n}$
                                                    \\ \hline
  $C$ & $Sp_{2k}(F)$ & $M_{k,n}$ & $(g,a) \mapsto ga$ &
  $so_{2n}$ \\ \hline
  \end{tabular}
\end{center}
\vspace{0.5cm} In case A we set $n=p+q$.  If $p=0$ (resp. $q=0$)
we set $V = M_{k,n}$ (resp. $V = M_{n,k}$). In all cases the
algebra of invariant differential operators $\mcD(V)^K$ is
generated (as an associative algebra) by a Lie algebra isomorphic
to $\fg$. Hence there is a surjective homomorphism
\begin{equation}\label{equation_phi}
    \phi: U(\fg) \longrightarrow \mcD(V)^K.
\end{equation}
The image of $\fg$ under $\phi$ is described explicitly on pages
69--70 of \cite{LS}.

Let $T$ be a maximal torus of $K$.  In \cite{MR}, it is shown that
if zero is not a $T$-weight of $V$ then $\mcD(V)^K$ has enough
finite dimensional representations.  Thus if $\fg$ is semisimple,
then $\mcD(V)^K$ is FCR since it is an image of $\mcU(\fg)$.  We
remark that zero is not a weight of $V$ in cases A and C of the
above table and in case B for even $k$.

There is a multiplicity free decomposition of $\mcO(V)$ as a $F[K]
\otimes \mcU(\fg)$-module (\cite{GW} Theorem 4.5.14). Furthermore,
as a $\mcU(\fg)$--module, $\mcO(V)$ is a direct sum of simple
highest weight modules. Let
\[ \La = \{ \la \in \fh^*| \mcO(V) \mbox{ has a $\fg$-submodule isomorphic to $L(\la)$} \}.
\]
If $\Omega$ is a radical ideal of $\mcS$, and $\Omega_1, \cdots,
\Omega_t$ are the prime ideals of $\mcS$ that are minimal over
$\Omega$, we set $I_\Omega = \cap_{i=1}^t I_{\Omega_i}$.

\begin{lem}\label{lem_A}
If $\Omega$ is the radical ideal of $\mcS$ such that $\mcV(\Omega) =
\overline{\La}$, the Zariski closure of $\La$, then
\begin{equation}
    \ker \phi = I_{\Omega}.
\end{equation}
\end{lem}
\begin{proof}
Suppose that $\overline{\La} = X_1 \cup \cdots \cup X_t$ is the
decomposition of $\overline{\La}$ into irreducible components, and
for $1 \leq i \leq t$, let $\Omega_i$ be the prime ideals of $\mcS$
such that $\mcV(\Omega_i) = X_i$.  By Lemma \ref{lem_density},
$\La_i = \La \cap X_i$ is dense in $X_i$, and so by (\cite{WS}
Proposition 1) we have:
\[
    I_{\Omega_i} = \bigcap_{\la \in \La_i} ann_{\mcU} L(\la)
\]
for $1 \leq i \leq t$.  Therefore, since $\mcO(V)$ is a faithful
$\mcD(V)^K$-module and $\La = \bigcup_{i=1}^t \La_i$ we have:
\begin{eqnarray*}
\ker \phi &=& \{u \in \mcU| \phi(u) \mcO(V) = 0 \} \\
          &=& \bigcap_{\la \in \La} ann_{\mcU} L(\la) \\
          &=& \bigcap_{i=1}^t \bigcap_{\la \in \La_i} ann_{\mcU}
                L(\la) \\
          &=& I_{\Omega}.
\end{eqnarray*}
\end{proof}

Note that $I_\Omega$ is a completely prime ideal.  We identify a
situation where it is zero.

\begin{thm}\label{thm_last}
Assume that rank $\fg \leq$ rank $K$ then $\mcU \cong \mcD(V)^K$.
\end{thm}
\begin{proof}
Let $\{\mcU_N\}$ be the usual filtration on $\mcU$ and $\{D_N\}$
the Bernstein filtration on $\mcD(V)$.  Then $K$ preserves each
$D_N$, so acts on $gr \mcD(V)$ and $gr \left( \mcD(V)^K \right) =
\left( gr \mcD(V) \right)^K$.  Now $\fg$ maps onto $D_2^K$ and
this induces a surjection $\mcU_N \rightarrow D_{2N}^K$ with
kernel $I_\Omega \cap \mcU_N$.  Passing to the associated graded
rings we obtain a surjection
\[\mcS(\fg) = gr \mcU \rightarrow gr
\mcD(V)^K=\mcS(V \oplus V^*)^K\] with kernel $gr I_\Omega$.  Then
if rank $K \geq $ rank $\fg$ then we can apply the Second
Fundamental Theorem of Invariant Theory in the free case
(\cite{GW} Corollary 4.2.5) to $V \oplus V^*$.  We conclude that
$gr I_\Omega = (0)$ so $I_\Omega=(0)$.
\end{proof}

Our next aim is to describe the irreducible constituents of
$\mcO(V)$ as a $K$-module and as a $\mcU(\f g)$-module.  With this
goal in mind we set up some standard notation.

By a partition we will mean a finite sequence of weakly decreasing
non-negative integers. We will denote partitions by lower case Greek
letters.  The \emph{length} of a partition is denoted $\ell(\a) =
\max \{i| \a_i > 0 \}$.  The \emph{conjugate} of a partition $\la$
is a partition $\la^\prime$ whose $i^{th}$ part is given by
$|\{j|\la_j \geq i \} |$.  Note that $\ell(\la) = \la^\prime_1$.

The highest weights of representations of $\fg$ (resp. $K$) are
given by $m$-tuples where $m$ is the rank of $\fg$ (resp. $K$).  The
$m$-tuple $(a_1, \cdots, a_m)$ corresponds to the weight
$\sum_{i=1}^m a_i \epsilon_i$ where the $\epsilon_i$ are defined in
Section 2.3.1 of \cite{GW}.

At this point we consider the cases separately.

\section{Case A}\label{sec_CASE_A}
As a basis for $\fh$ we take $\{E_1, \cdots, E_n\}$ where $E_i$ is
the $n \times n$ matrix with a 1 in the row $i$, column $i$ entry
and zeros elsewhere.

Given non-negative integer partitions $\a$ and $\b$ with $\ell(\a)
\leq p$, $\ell(\b) \leq q$, and $\ell(\a) + \ell(\b) \leq k$, let
$V_{(\a,\b)}$ denote the irreducible $\fgl_n(F)$ module with highest
weight:
\[
\begin{array}{ccc}
    (\a,\b)_{\fg} &:=& \underbrace{(-k-\a_p,-k-\a_{p-1},\cdots,-k-\a_1, \b_1,\b_2,\cdots, \b_q)} .\\
                   & & n
\end{array}
\]

Let $V^{(\a,\b)}$ denote the irreducible representation of $GL_k$
with highest weight:
\[
\begin{array}{ccc}
    (\a,\b)_K &:=& \underbrace{(\a_1, \a_2, \cdots, \a_p,0,\cdots,0,-\b_q, \cdots, -\b_2, -\b_1)}. \\
              &  & k
\end{array}
\]

\begin{lem}
There is a multiplicity free decomposition under the joint action
of $K$ and $\fg$ we have the multiplicity free decomposition:
\[
    \mcO(V) = \bigoplus_{} V^{(\a,\b)} \otimes
    V_{(\a,\b)},
\]
where the direct sum is over the set
\[W_k(p,q) := \left\{ (\a, \b) \left|
\begin{array}{lll}
\hbox{$\a$ and $\b$ partitions such that} \\
{\ell(\a) \leq p, \, \ell(\b) \leq q, \quad \hbox{and}} \\
{\ell(\a)+\ell(\b) \leq k} \\
\end{array}
\right. \right\} \]
\end{lem}
\begin{proof} See \cite{Annals} p. 356.
\end{proof}

Define $\La_k(p,q) = \{ (\a,\b)_\fg | (\a,\b) \in W_k(p,q) \}$,
and let $Z_k(p,q) := \overline{\La_k(p,q)}$ denote the Zariski
closure of these weights.

If $k < n$ define:
\[
\Omega_m := \left\{%
\begin{array}{ll}
    (E_{k+1}, E_{k+2}, \cdots, E_n), & \hbox{if $m=0$;}                                                    \\
    (E_1+k, \cdots, E_m +k, E_{m+k+1}, \cdots, E_n), & \hbox{if
    $0<m<n-k$;} \\
    (E_1 + k, \cdots, E_{n-k}+k), & \hbox{if $m=n-k$.}                                                         \\
\end{array}%
\right.
\]


\begin{lem}

The decomposition of $Z_k(p,q)$ into irreducible components is as
follows:

\begin{enumerate}
\item $Z_k(p,q) = \fh^*$, if $n\leq k$. \item If $k < n \leq 2k$,
then
\[
    Z_k(p,q) = \bigcup_{m \in \Phi_p} \mcV(\Omega_m)
\]
where:
\[
\Phi_p = \left\{%
\begin{array}{ll}
    \{  0, \cdots,   p \}, & \hbox{if $p \leq n-k$;}        \\
    \{  0, \cdots, n-k \}, & \hbox{if $n-k \leq p \leq k$;} \\
    \{p-k, \cdots, n-k \}, & \hbox{if $k \leq p \leq n$.}   \\
\end{array}%
\right.
\]

\item If $2k \leq n$ then
\[Z_k(p,q) = \bigcup_{m \in \Phi_p} \mcV(\Omega_m) \]
where:
\[ \Phi_p =
\left\{%
\begin{array}{ll}
    \{  0, \cdots,   p \}, & \hbox{if $p \leq k$;} \\
    \{p-k, \cdots,   p \}, & \hbox{if $k \leq p \leq n-k$;} \\
    \{p-k, \cdots, n-k \}, & \hbox{if $p \geq n-k$.}
\end{array} \right.\]
\end{enumerate}
\end{lem}
\begin{proof}
Statement (1) follows because if $n \leq k$, then $\mcO(V)$ contains
the irreducible $\fgl_n(F)$ module with highest weight
$(\a,\b)_{\fg}$ for any partitions $\a,\b$ with $\ell(\a) \leq p$,
and $\ell(\b) \leq q.$ For statements (2) and (3) we define, for $0
\leq i \leq p$
\[ W_k^i(p,q) = \{ (\a,\b) \in W_k(p,q)| \ell(\a) = p-i \},\]
and
\[ \La_k^i(p,q) = \{ (\a,\b)_\fg | (\a,\b) \in W_k(p,q) \}. \]

Suppose $(\a,\b) \in W^i_k(p,q)$.  Then, $0 \leq i \leq p$ since $0
\leq \ell(\a) \leq p$.  Also, $\ell(\a) \leq k$ implies that $p-k
\leq i$.  In addition, the set of weights $(\a,\b)_\fg$ such that
$\ell(\a) = p-i$ and $\ell(\b) < k-\ell(\a) = k+i-p$ is contained in
the Zariski closure of the set of weights $(\a,\b)_\fg$ with
$\ell(\b) = k+i - p.$  Since we require that $\ell(\b) \leq q$, we
can restrict our attention to the highest weights $(\a,\b)_\fg$ such
that $\ell(\a) = p-i$ and $k+i-p \leq q$.  That is, $i \leq n-k.$ It
follows that $Z_k(p,q) = \bigcup \overline{\La^i_k(p,q)}$ where the
union is over all $i$ such that $ \max(0,p-k) \leq i \leq \min(p,
n-k)$.  It is easy to check that for these values of $i$,
$\overline{\La^i_k(p,q)} = \mcV(\Omega_i)$. The result follows by
considering individual cases.
\end{proof}

\begin{prop}
If $n$ and $k$ are fixed, then for $0 \leq i < j \leq n-k$ we have
$I_{\Omega_i} = I_{\Omega_j}$. Moreover, $\Omega_0$ and
$\Omega_{n-k}$ are strongly dominant.
\end{prop}
\begin{proof}
To see that $\Omega_0$ is strongly dominant note that
\[ \sum_{i = 1}^k \bbN \varp_i \subseteq \mcV(\Omega_0) = \sum_{i = 1}^k  \rF \varp_i. \]
Similarly $\Omega_{n-k}$ is strongly dominant. Now set $\Omega =
\Omega_{n-k}$, so that:
\[ \Omega = (E_1+k, E_2+k, \cdots, E_{n-k}+k), \]
and let $\la = \la_\Omega$.  For $0 \leq i \leq n-k$ we define
$w_i \in W = S_n$ by:
\[
    w_i(j) = \left\{%
\begin{array}{ll}
    j,       & \hbox{for $1     \leq j \leq   i$;} \\
    j+k,     & \hbox{for $i+1   \leq j \leq n-k$;} \\
    j-n+k+i, & \hbox{for $n-k+1 \leq j \leq   n$.} \\
\end{array}%
\right.
\]
Then, it is easy to see that $w_i.\Omega = \Omega_i$, and that $w_i
\in W^\la$. The result follows from the Lemma \ref{lem_c->b}.
\end{proof}

\begin{rem}
    Explicit generators for $I_{\Omega_0}$ are given in \cite{VP}
    using Capelli identities.
\end{rem}

\section{Cases B and C}\label{sec_CASE_BC}

We first turn our attention to Case B.

\begin{lem}
There is a multiplicity free decomposition under the joint action
of $K$ and $\fg$ given by
\[
    \mcO(V) = \bigoplus_{\mu \in M} V^{\mu} \otimes
    V_{\mu},
\]
where we index the summands by the partitions $\mu$ in the set:
\[M := \left\{
\hbox{$\mu$ is a partition such that} \; \mu'_1 +\mu'_2 \leq k \;
\mbox{and} \; \ell := \mu_1^\prime \leq n \right \} .\] The highest
weight of the $sp_{2n}$-module, $V_\mu$ is given by:
\[
\begin{array}{ccc}
 & \underbrace{(-\frac{k}{2},-\frac{k}{2},\cdots,-\frac{k}{2}-\mu_{\ell},\cdots,-\frac{k}{2}-\mu_1)} \\
                   & n
\end{array}
\]
\end{lem}
\begin{proof}
See \cite{Annals} p. 353.
\end{proof}

We describe the Zariski closure, $Z$ in $\fh ^*$ of the set of
weights of $\mcO (V)$ as a $\fg$-module.  If $k \geq 2n$ then $Z =
\fh^*$, so we consider the case where $k < 2n$.  For $k$ odd, we can
see using Theorem \ref{thm_ideal-containment} that $\mcD(V)^K$ has
no finite dimensional modules, so we consider the $k$ even case.

Suppose that $k = 2p$.  Fix $i$ with $1 \leq i \leq \min(p,n-p)$
and set $r = n-p-i$, $s = n-p+i.$ Then define
\[ \Omega_i  =  (H_1, H_2, \cdots, H_{r-1}, H_r +1, H_{r+1}, \cdots,
H_{s-1},H_s + H_{s+1} + \cdots H_{n} +p +1 ).\] We also set
\[ \Omega_0   = (H_1,\cdots,
H_{n -p -1}, H_{n-p} + H_{n -p +1}+ \cdots + H_{n} +p).\]

\begin{lem}\label{lem_decomposition}
The decomposition of $Z$ into irreducible components is given by
\[ Z = \bigcup_{i = 0}^{\min(p,n-p)} \mcV(\Omega_i). \]
\end{lem}
\begin{proof}  For $0 \leq i \leq \min(p,n-p)$, set
\[ M_i = \{ \mu \mbox{ is a partition such that } \mu_1^\prime + \mu_2^\prime \leq k, \mbox{ and } \mu_1^\prime = p+ i  \}. \]
and let $\La_i \subseteq \fh^*$ be the set of highest weights of the
modules $V_\mu$ with $\mu \in M_i$.  Note that if $\mu_1^\prime +
\mu_2^\prime \leq k=2p$ and $\mu_1^\prime < p$ then the highest
weight of $V_\mu$ is in $\overline \La_0$.

If $0 \leq i < j \leq \min(p,n-p)$ neither of the ideals $\Omega_i,
\Omega_j$ is contained in the other.  Therefore, it suffices to show
that the Zariski closure of $\La_i$ is $\mcV(\Omega_i)$.  Suppose
that $\mu \in M_i$ and let $\a = \sum_{i=1}^n a_i \epsilon_i$ be the
highest weight of $V_\mu$.  Then $a_i = -p$ for $1 \leq i \leq r =
n-p-i$, and $a_i = -p - 1$ for $r+1 \leq i \leq s = n-p+i$.  The
result follows easily.
\end{proof}

Note that $\mbox{GK-dim} \; \mcS /\Omega_{i} = p - i$, where
GK-dim is the Gelfand-Kirillov dimension.  Let $I_i =
I_{\Omega_{i}}$ for $0 \leq i \leq p.$

Recall the homomorphism $\phi$ from equation (\ref{equation_phi}).
\begin{lem}\label{lem_B}
We have $\ker \phi = I_0.$
\end{lem}
\begin{proof}
By Lemma \ref{lem_A}, $\ker \phi = \cap_{i=0}^p I_i$.  Since $\ker
\phi$ is prime, it follows that $I_i \subseteq \ker \phi \subseteq
I_0$ for some $i$.  Hence $I_i \cap \mcZ \subseteq I_0 \cap \mcZ$,
that is $\Omega_i \cap \mcZ \subseteq \Omega_0 \cap \mcZ$.  If
$i\geq 1$ then since $\mcS$ is a finitely generated $\mcZ$-module we
would have:
\begin{eqnarray*}
p-i &=   & \mbox{GK-dim}(S/\Omega_i) = \mbox{GK-dim}(Z/{Z\cap \Omega_i}) \\
    &\geq& \mbox{GK-dim}(Z/{Z\cap \Omega_0}) = p
\end{eqnarray*}
a contradiction.
\end{proof}
We show next that $\mcU/I_0$ is FCR.  To do this define \[\Omega =
(H_{p +1}, \cdots, H_{n}),\] and define $w \in W$ so that if $a =
(a_{1}, \cdots, a_{n})\in \fh^*$, we have $w(a) = b,$ where:

\[
b_i =
\left\{%
\begin{array}{ll}
    a_{i+p  },  & \hbox{for $1     \leq i \leq n-p$;} \\
    a_{i+p-n},  & \hbox{for $n-p+1 \leq i \leq n  $.} \\
\end{array}%
\right.
\]

 It is easy to check that
$\Omega$ is strongly dominant, $w.\Omega = \Omega_{0}$ and $w \in
W^{\la}$ where $\la = \la_{\Omega}.$  Thus by the main theorem, $I_0
= I_{\Omega}$ and $\mcU_{\Omega}$ is FCR.  We expect that
\[ I_0 \subset I_1 \subset \cdots \subset
I_p\] with all inclusions strict.  If $p=1$ this follows from
Lemmas \ref{lem_A} and \ref{lem_B}.



Case C is relatively uninteresting, the Zariski closure
$\mcV(\Omega)$ of the set of highest weight occurring in $\mcO(V)$
as a $\fg$-module is irreducible.  Also $\mcU_\Omega \cong D(V)^K$
has enough finite dimensional modules from \cite{MR} and thus is
FCR.  Furthermore, $\mcV(\Omega) = \fh^*$ if and only if $n \leq
k$.

\bibliographystyle{alpha}

\begin{thebibliography}{How89}

\bibitem[BJ01]{BJ}
Walter Borho and Anthony Joseph.
\newblock {S}heets and topology of primitive spectra for
  semisimiple {L}ie algebras {\em J. Algebra}, 259(1):310--311,
  2003.;
\newblock Corrigendum to   {J}. {A}lgebra {\bf 244} (2001), no. 1,
  76--167.


\bibitem[Dix96]{JD}
J. Dixmier. {\em Enveloping algebras.} Revised reprint of the 1977
translation. American Mathematical Society, Providence, RI, 1996.


\bibitem[EW04]{Annals}
T.~J. Enright and J.~F. Willenbring.
\newblock Hilbert series, {H}owe duality and branching for classical groups.
\newblock {\em Ann. of Math. (2)}, 159(1):337--375, 2004.

\bibitem[GW98]{GW}
R.~Goodman and N.R. Wallach.
\newblock {\em Representations and invariants of the classical groups}.
\newblock Cambridge University Press, Cambridge, 1998.



\bibitem[How89]{RH}
R.~Howe.
\newblock Remarks on classical invariant theory.
\newblock {\em Trans. Amer. Math. Soc.}, 313(2):539--570, 1989.

\bibitem[Hu90]{Hu}
J.E.~Humphreys.
\newblock {\it Reflection groups and Coxeter groups}.
\newblock Cambridge Univ. Press, Cambridge, 1990.

\bibitem[Jan83]{Ja}
J.~C. Jantzen.
\newblock {\em Einh\"ullende {A}lgebren halbeinfacher {L}ie-{A}lgebren},
  volume~3 of {\em Ergebnisse der Mathematik und ihrer Grenzgebiete (3)
  [Results in Mathematics and Related Areas (3)]}.
\newblock Springer-Verlag, Berlin, 1983.


\bibitem[Jo77]{Jo77}
A. Joseph. {A characteristic variety for the primitive spectrum of a
semisimple Lie algebra.} {\em Non-commutative harmonic analysis}
(Actes Colloq., Marseille-Luminy, 1976), pp. 102--118. Lecture Notes
in Math., Vol. 587, Springer, Berlin, 1977.


\bibitem[KS94]{KS}
H.~Kraft and L.~W. Small.
\newblock Invariant algebras and completely reducible representations.
\newblock {\em Math. Res. Lett.}, 1(3):297--307, 1994.

\bibitem[LS89]{LS}
T.~Levasseur and J.~T. Stafford.
\newblock Rings of differential operators on classical rings of invariants.
\newblock {\em Mem. Amer. Math. Soc.}, 81(412):vi+117, 1989.

\bibitem[MR]{MR}
I.~M. Musson and S.~L. Rueda.
\newblock {Finite dimensional representations of invariant differential
  operators}.
\newblock {Trans. Amer. Math. Soc. 357 (2004), no. 7, 2739-2752.}

\bibitem[P05]{P05}
Nikolaos Papalexiou.
\newblock{ Order relations on the induced and prime ideals in the enveloping algebra of a semisimple Lie algebra}
\newblock{ Journal of Algebra 302 (2006) 1--15.}

\bibitem[P]{P}
Nikolaos Papalexiou.
\newblock{ On the prime spectrum of the enveloping algebra and characteristic varieties.}
\newblock{ To Appear in the Journal of Algebra and its Applications.}

\bibitem[Pro04]{VP}
V. Protsak.
\newblock{Rank ideals and Capelli identities}.
\newblock {J. Algebra 273 (2004), no. 2, 686--699.}

\bibitem[Soe90]{WS}
W. Soergel.
\newblock The prime spectrum of the enveloping algebra of a reductive {L}ie
  algebra.
\newblock {\em Math. Z.}, 204(4):559--581, 1990.

\bibitem[Soe]{WS-web}
W. Soergel. \\
\newblock \url{http://home.mathematik.uni-freiburg.de/soergel/PReprints/Korrektur.pdf}

\end{thebibliography}

\def\cprime{$'$} \def\cprime{$'$}

\end{document}